\documentclass[11pt]{amsart}
\usepackage{amssymb, amstext, amscd, amsmath}
\usepackage[all]{xy}
\theoremstyle{plain}
\newtheorem{theorem}{Theorem}[section]
\newtheorem{corollary}[theorem]{Corollary}
\newtheorem{proposition}[theorem]{Proposition}
\newtheorem{lemma}[theorem]{Lemma}

\theoremstyle{definition}
\newtheorem{definition}[theorem]{Definition}
\newtheorem{definitions}[theorem]{Definitions}

\newtheorem{remark}[theorem]{Remark}

\makeatletter
\def\@cite#1#2{{\m@th\upshape\bfseries%
[{#1\if@tempswa{\m@th\upshape\mdseries, #2}\fi}]}} \makeatother

\newcommand{\bbC}{{\mathbb{C}}}

\newcommand{\bbN}{{\mathbb{N}}}

\newcommand{\B}{{\mathcal{B}}}

\newcommand{\I}{{\mathcal{I}}}
\newcommand{\J}{{\mathcal{J}}}

\newcommand{\fI}{{\mathfrak{I}}}

\renewcommand{\phi}{\varphi}
\newcommand{\upchi}{{\raise.35ex\hbox{\ensuremath{\chi}}}}
\def\gs{\sigma}

\def\gl{\lambda}

\def\eps{\varepsilon}

\newcommand{\foral}{\text{ for all }}

\newcommand{\id}{{\operatorname{id}}}

\newcommand{\ext}{\operatorname{ext}}

\newcommand\Span{\mathop{\rm span}}

\newcommand\diag{\mathop{\rm diag}}

\newcommand{\ca}{\mathrm{C}^*}
\newcommand{\cenv}{\mathrm{C}^*_{\text{env}}}

\newcommand{\sca}[1]{\left\langle#1\right\rangle}
\newcommand{\nor}[1]{\left\Vert #1\right\Vert}
\begin{document}

\title{The \v{S}ilov Boundary for operator spaces}

\author{Evgenios T.A. Kakariadis}
\address{Pure Mathematics Department, University of Waterloo,
   On N2L3G1, Canada}
\email{ekakaria@uwaterloo.ca}

\thanks{The author was partially supported by the Fields Institute for Research in the Mathematical Sciences.}

\maketitle

%%%%%%%%%%%%%%%%%%%%%%
\begin{abstract}
Motivated by the recent interest in the examination of unital completely positive maps and their effects in $\ca$-theory, we revisit an older result concerning the existence of the \v{S}ilov ideal. The direct proof of Hamana's theorem for the existence of an injective envelope for a unital operator subspace $X$ of some $\B(H)$ that we provide implies that the \v{S}ilov ideal is the intersection of $\ca(X)$ with any maximal boundary operator subsystem in $\B(H)$. As an immediate consequence we deduce that the \v{S}ilov ideal is the biggest boundary \emph{operator subsystem} for $X$ in $\ca(X)$.

The new proof of the existence of the \v{S}ilov ideal that we give does not use the existence of maximal dilations, provided by Drits-\break chel and McCullough, and so it is independent of the one given by Arveson. As a consequence, the \v{S}ilov ideal can be seen as the set that contains \emph{the abnormalities} in a $\ca$-cover $(C,\iota)$ of $X$ for all the extensions of the identity map $\id_{\iota(X)}$. The interpretation of our results in terms of ucp maps characterizes the maximal boundary subsystems of $X$ in $\B(H)$ as kernels of $X$-projections that induce completely minimal $X$-seminorms; equivalently, $X$-minimal projections with range being an injective envelope, that we view from now on as the \v{S}ilov boundary for $X$.
\end{abstract}

%%%%%%%%%%%%%%%%%%%%%%
\section{Introduction}
%%%%%%%%%%%%%%%%%%%%%%

Let $X,Y$ be linear spaces and $\phi\colon X \rightarrow Y$ a linear map. We define $\phi_\nu: =\id_\nu\otimes\phi\colon M_\nu(X)\rightarrow \B(H^\nu)$ by $\phi_\nu([a_{ij}])=[\phi(a_{ij})]$.

An (abstract) operator space is a pair $(X, \{\|\cdot\|_{\nu}\}_{\nu\geq 1})$, consisting of a vector space, and a norm on $M_\nu(X)$ for all $\nu\in\bbN$, such that there exists a linear map $u\colon X\rightarrow \B(H)$ (where $H$ is a Hilbert space) such that every $u_\nu$ is an isometry; equivalently, Ruan's axioms hold for the sequence of norms. In this case we call the sequence $\{\|\cdot\|_{\nu}\}_{\nu\geq 1}$ an operator space structure on the vector space $X$. Throughout this paper $X$ is assumed unital, i.e., there is an element $e\in X$ such that $u(e)=I_H$, and the morphisms will always be unital.

If $X$ is a linear subspace of a $C^*$-algebra $C$, then $X$ is an operator space with the matrix norm structure inherited by a faithful representation of $C$. By Ruan's theorem we can always assume that an operator space sits inside a $\ca$-algebra.

An \emph{operator system} is a selfadjoint linear subspace $S$ of a unital $C^*$-algebra, that contains the unit. (There is an alternative way of defining abstract operator systems, which we will not make use of.)

Let $X$ be a unital operator space and $\phi\colon X\rightarrow \B(H)$ be a linear map.  We define $\phi_\nu: =\id_\nu\otimes\phi\colon M_\nu(X)\rightarrow \B(H^\nu)$ by $\phi_\nu([a_{ij}])=[\phi(a_{ij})]$. We call $\phi$ \emph{unital completely positive $($ucp$)$, unital completely contractive $($ucc$)$ or unital completely isometric $($ucis$)$} if $\phi_\nu$ is respectively positive, contractive or an isometry, for every $\nu\in\bbN$.

We can use the decomposition of an element $x\in S_{sa}$ in two positive elements in $S$, i.e., $x= (1\|x\|+x)/2 - (1\|x\|-x)/2$, to prove that a linear map $\phi\colon S\rightarrow \B(H)$ is ucp if, and only if, it is ucc. Also, if $X$ is a unital operator space sitting in a $\ca$-algebra and $\phi\colon X \rightarrow \B(H)$ is a ucc map, then $X+X^*$ is an operator system and there is a unique ucp extension $\tilde{\phi}\colon X+X^* \rightarrow \B(H)$ of $\phi$.

Given two ucc maps $\phi_k\colon X\rightarrow \B(H_k)$, $k=1,2$, we write $\phi_1\leq\phi_2$ if $H_1\subseteq H_2$ and $P_{H_1}\phi_2(x)|_{H_1}=\phi_1(x)$, $x\in X$; $\phi_2$ is called a \emph{dilation} of $\phi_1$ and $\phi_1$ is called a \emph{compression} of $\phi_2$. In general, we can have the following scheme that gives the existence of dilations of a ucc map. Let $X\subseteq\B(H)$ be a unital operator space and $\phi\colon X\rightarrow \B(K)$ a ucc map. \emph{Arveson's Extension Theorem} implies that there is a ucp (thus ucc) map $\psi\colon \B(H) \rightarrow \B(K)$ extending $\phi$. By applying \emph{Stinespring's Dilation Theorem} on $\psi$, there is a Hilbert space $W\supseteq K$ and a unital representation $\pi\colon \B(H) \rightarrow\B(W)$ such that $\psi(c)=P_K\pi(c)|_K$, for every $c\in \B(H)$. Hence, $\pi|_X$ is a dilation of $\phi$.

Given an operator space $X$, a natural question to ask is which is the \emph{smallest} (in some sense) $\ca$-algebra $C$ for which there is a ucis map $\phi\colon X \rightarrow C$, i.e., \emph{the $\ca$-envelope of $X$}. \emph{A $\ca$-cover} for $X$ is a pair $(C,\iota)$ where $\iota\colon X \rightarrow C$ is a ucis map amd $C=\ca(\iota(X))$.

%%%%%%%%%%%%%%%%%%%%%%
\begin{definition}
Let $X$ be a unital operator space. The $\ca$-envelope of $X$ is the $\ca$-cover $(\cenv(X),\iota)$  with the following (universal) property:
\emph{ for every $\ca$-cover $(C,j)$ there exists a unique $*$-epimorphism $\Phi\colon  C \rightarrow \cenv(X)$, such
that $\Phi(j(x))=\iota(x)$, for every $x\in X$}.
\end{definition}

The existence of the $\ca$-envelope was first proved by Arveson, in the case where there were enough boundary representations. The first proof for the general case was given by Hamana in \cite{Ham79}. Twenty five years later Dritschell and McCullough gave an independent proof in \cite{DrMc05} for the existence of the $\ca$-envelope. We should remark here that the original versions of these theorems were stated in terms of operator algebras or operator systems. A moment of clarity shows that one can easily reformulate these theorems for (unital) operator spaces $X$, simply by mimicking the simplified proof in \cite{Ar06}.

The key step of proving the existence of the $\ca$-envelope in \cite{DrMc05} was the proof of the existence of \emph{a maximal} representation for $X$. The following definitions are equivalent.

%%%%%%%%%%%%%%%%%%%%%%
\begin{definitions}
(i) A ucc map $\phi\colon X\rightarrow \B(H)$ is said to be \emph{maximal} if it has no nontrivial dilations, i.e. $\phi'\geq \phi \Rightarrow \phi'=\phi\oplus\psi$, for some ucc map $\psi$,

\noindent (ii) A ucc map $\pi\colon X\rightarrow \B(H)$ is said to have the \emph{unique extension property} if

1. $\pi$ has a unique completely positive extension $\tilde{\pi}\colon C^*(X)\rightarrow\B(H)$, and

2. $\tilde{\pi}\colon C^*(X)\rightarrow \B(H)$ is a representation of $C^*(X)$ on $H$.

\noindent(iii) A ucc map $\phi\colon X \rightarrow \B(H)$ is called a \emph{$\partial$-representation} if for any dilation $\nu \geq \phi$ the Hilbert space $H$ is $\nu(X)$-reducing.
\end{definitions}

%%%%%%%%%%%%%%%%%%%%%%
\begin{theorem} \textup{(Dritschel-McCullough)}
Let $X$ be an operator space in a $\B(H)$. Then the identity map $\id\colon X \rightarrow \B(H)$ has a dilation $\nu \colon X \rightarrow \B(K)$ that is maximal. Hence, $\ca(\nu(X)) \simeq \cenv(X)$.
\end{theorem}

In contrast, the direction for the proof in \cite{Ham79} is completely different and has an algebraic flavor. It is based on proving the existence of \emph{an injective envelope} by using the notion of minimal $X$-seminorms.

%%%%%%%%%%%%%%%%%%%%%%
\begin{definition}
An operator space $E$ is called \emph{injective} if for any pair $Z,Y$ of operator spaces such that $Z\subseteq Y$, and every ucc map $\phi\colon Z \rightarrow E$, there exists a ucc map $\psi\colon Y \rightarrow E$ that extends $\phi$.
\end{definition}

If $\iota \colon X \rightarrow E$ is a ucis in an operator space $E$, then the pair $(E,\iota )$ is called an extension of $X$.  We say that an extension $(E,\iota )$ is \emph{rigid} if $\id_E$ is the only ucc map $E \rightarrow E$ that extends the identity map on $\iota (X)$. We say that an extension $(E,\iota )$ is \emph{essential} if whenever $\phi\colon E \rightarrow Z$ is a ucc map into another operator space $Z$ such that $\phi\circ \iota $ is a complete isometry, then $\phi$ is a complete isometry. We say that $(E,\iota)$ is an injective \emph{envelope} of $X$ if $E$ is injective and there is no injective subspace of $E$ containing $\iota(X)$. One can prove that an injective extension $(E,\iota)$ is an envelope of $X$ if and only if it is rigid if and only if it is essential. Also, if $(E,\iota)$ is an injective envelope and $\phi\colon E \rightarrow \B(K)$ is a ucc map such that the restriction of $\phi$ to $\iota(X)$ is ucis, then $\phi$ is a ucis map and $(\phi(E), \phi\circ \iota)$ is an injective envelope for $X$ in $\B(K)$.

%%%%%%%%%%%%%%%%%%%%%%
\begin{theorem}
\textup{(Hamana)} Let $X$ be an operator space in a $\B(H)$. Then there is an injective envelope $(E,\iota)$ of $X$. Thus $\ca(\iota(X))\simeq \cenv(X)$.
\end{theorem}

We pinpoint two lemmas concerning injective envelopes of Hamana's theory that we are going to use in the following sections.

%%%%%%%%%%%%%%%%%%%%%%
\begin{lemma}\label{nec un in en}
Let $(S_k,\iota_k)$, $k=1,2$, be injective envelopes for an operator space $X$. Then the mapping $\iota_1(x)\rightarrow \iota_2(x)$ extends to a ucc map $\phi\colon S_1 \rightarrow S_2$ which is a necessarily unique ucis onto map.
\end{lemma}
\begin{proof} By injectivity of $S_2$ there is a ucc extension $\phi\colon S_1\rightarrow S_2$ that fixes $X$ elementwise, i.e., $\phi(\iota_1(x))=\iota_2(x)$, for all $x\in X$. Since $S_1$ is also an essential envelope $\phi$ is a ucis map and $\phi(S_1)$ is an injective envelope in $S_2$. Therefore $\phi$ is also onto $S_2$. Now let $\psi\colon S_1 \rightarrow S_2$ be a ucis map such that $\psi(\iota_1(x))=\iota_2(x)$, for all $x\in X$. Then $\psi$ is also a ucis map onto $S_2$. Moreover the restriction of $\psi^{-1}\circ \phi$ to $\iota_1(X)$ is the identity mapping, hence $\psi^{-1}\circ \phi=\id$ on $S_1$, by rigidity of $S_1$. Therefore $\psi=\phi$.
\end{proof}

%%%%%%%%%%%%%%%%%%%%%%
\begin{lemma}\label{nec onto in en}
Let $(E,j)$ be an injective extension and $(S,\iota)$ be an injective envelope for an operator space $X\subseteq \B(H)$. Then any ucc map $\phi\colon E \rightarrow S$ such that $\phi(j(x))=\iota(x)$, for all $x\in X$, is onto $S$. Moreover, if $E=\B(H)$ and $(S,\id_X)$ is the injective envelope, then $\phi$ is a projection.
\end{lemma}
\begin{proof} Let such a map $\phi\colon E \rightarrow S$ and $\gs\colon S \rightarrow E$ be an extension of the mapping $\iota(x)\mapsto j(x)$, for all $x\in X$. Then $\phi\circ \gs(\iota(x))=\iota(x)$ for all $x\in X$ and the range of $\phi\circ \gs$ is in $S$. Thus $\phi\circ \gs=\id$ by rigidity of $S$, and the proof of the first statement is complete.

If $S\subseteq \B(H)$, then the restriction of $\phi\colon S\equiv \phi(\B(H)) \rightarrow S$ to $X=\id(X)=\phi(X)$ is the identity map, therefore $\phi^2(a)=\phi(\phi(a))=\id(\phi(a))=\phi(a)$, for all $a\in \B(H)$, by rigidity of $S$.
\end{proof}

The approaches of \cite{Ham79} and \cite{DrMc05} gave independently the existence of the $\ca$-envelope, thus the existence of a second object, the \emph{\v{S}ilov ideal}. As a result the \v{S}ilov ideal is described as the kernel of a necessarily unique $*$-epimorphism.

%%%%%%%%%%%%%%%%%%%%%%
\begin{definition}
If $\iota\colon X \rightarrow C$ is a ucis map and $C=\ca(\iota(X))$, then an ideal $I$ of $\ca(\iota(X))$ is called boundary if the restriction of the natural $*$-epimorphism $q_I\colon C \rightarrow C/I$ to $\iota(X)$ is a ucis map. The biggest boundary ideal in $\ca(\iota(X))$ is called \emph{the \v{S}ilov ideal of $\iota(X)$ in $C$}.
\end{definition}

It appears that the \v{S}ilov ideal is a very tractable tool for finding the $\ca$-envelope in recent papers (and I will avoid making any advertisement here, as it is irrelevant to our subject). The crucial remark used in some of these cases is that the $\ca$-envelope contains no non-trivial boundary ideals. Indeed, if $I$ is the \v{S}ilov ideal in a $\ca$-cover $(C,\iota)$, then $C/I \simeq \cenv(X)$. Let us show here how the existence of the $\ca$-envelope implies the existence of the \v{S}ilov ideal. Note that by definition the \v{S}ilov ideal is unique.

%%%%%%%%%%%%%%%%%%%%%%
\begin{proposition}
If there exists a $\ca$-cover for an operator space $X$ that has the universal property of the $\ca$-envelope, then the \v{S}ilov ideal exists.
\end{proposition}
\begin{proof}
Assume that $(\cenv(X),\iota)$ has the universal property and let $(C,j)$ be a $\ca$-cover for $X$. Then there is a unique $*$-epimorphism $\Phi\colon C \rightarrow \cenv(X)$, such that $\Phi(j(x))=\iota(x)$ for all $x\in X$. Then, by the first theorem for $*$-isomorphisms $C/\ker\Phi \simeq \cenv(X)$, via the $*$-isomorphism
\[
\widehat\Phi(c+\ker\Phi)= \Phi(c), \foral c\in C.
\]
Since $\widehat\Phi$ is a $*$-isomorphism, hence a ucis map, the ideal $\ker\Phi$ is boundary.

Let $J$ be any boundary ideal in $C$. Then $(C/J, q_J\circ j)$, where $q_J\colon C \rightarrow C/J$ is the canonical $*$-epimorphism, is a $\ca$-cover for $X$. Therefore, by the universal property of $\cenv(X)$, there is a unique $*$-epimorphism $\Pi\colon C/J \rightarrow \cenv(X)$, such that $\Pi(q_J(j(x)))=\iota(x)$ for all $x\in X$. Then $\Pi\circ q_J(j(x)) = \Phi(j(x))$ for all $x\in X$, thus $\Pi\circ q_J=\Phi$, since $\Pi, q_J$ and $\Phi$ are $*$-homomorphisms and $C=\ca(j(X))$. Hence,
\[
\ker\Phi=\ker(\Pi\circ q_J) \supseteq \ker q_J= J.
\]
Therefore, the ideal $\ker\Phi$ contains all the boundary ideals in $C$. So it is the \v{S}ilov ideal.
\end{proof}

It is interesting that there is not a known proof of the converse without the additional use of the existence of maximal representations proved by Dritschell and McCullough \cite{DrMc05}, except from that given by Arveson in \cite{Ar06}, or the use of the existence of an injective envelope.  Yet, the advantage of a direct proof of the existence of the \v{S}ilov ideal provides additional information; the proof provided by Arveson in \cite{Ar06} characterizes the \v{S}ilov ideal as the kernel of (the unique extension of) a maximal representation. The proof, that we provide in what follows, gives additional characterizations of the \v{S}ilov ideal in terms of maximal boundary subsystems and/or kernels of minimal $X$-maps.

%%%%%%%%%%%%%%%%%%%%%%
\section{The proof}
%%%%%%%%%%%%%%%%%%%%%%

Let us fix a Hilbert space $H$ such that $X \subseteq \B(H)$ (completely isometrically). If $S\subseteq \B(H)$ is an operator system that contains $X$, we say that a selfadjoint subspace $\I$ of $S$ is \emph{a boundary subsystem} for $X$, if the restriction to $X$ of the quotient linear map $q_\I\colon S \rightarrow S/\I\subseteq \B(H)/\I$ is completely isometric. The matrix norms in $\B(H)/\I$ come from the identification of $M_\nu(\B(H)/\I)\simeq M_\nu(\B(H))/M_\nu(\I)$, i.e.,
\begin{align*}
\nor{[x_{ij}]+ M_\nu(\I)}_\nu=\nor{[x_{ij}+\I]}_\nu=\inf\{\nor{[x_{ij}+y_{ij}]}_\nu:y_{ij}\in \I\}
\end{align*}
It is trivial to see that boundary subsystems are never unital, since $X$ is unital. The translation of the invariance principle in our context is the following.

%%%%%%%%%%%%%%%%%%%%%%
\begin{proposition}\label{pull back proposition}
If $\I$ is a boundary subsystem for $X$ in a $\B(H)$ and $V/\I$ is a boundary subsystem of $q_\I(X)$ in $\B(H)/\I$, then $V$ is a boundary subsystem for $X$ in $\B(H)$.
\end{proposition}
\begin{proof} By contractivity of $q_\I$, for $x \in X$ and $v\in V$ we
get that
\begin{align*}
\nor{q_\I(x) + q_\I(v)}= \nor{q_\I(x +v)} \leq \nor{x +v}
\end{align*}
hence by taking the infimum over all $v\in V$
\begin{align*}
\nor{x}=\nor{q_\I(x)}&= \nor{q_\I(x) + V/\I} = \inf\{ \nor{q_\I(x) + q_\I(v)}\colon  v \in V\}\\
 & \leq \inf\{ \nor{x +v}\colon  v\in V\} =\nor{x + V} \leq \nor{x}.
\end{align*}
Thus $\nor{x}=\nor{x+V}$, for all $x\in X$, so the restriction of $q_V$ on $X$ is isometric. A similar argument for all matrix norms gives that $V$ is a boundary subsystem for $X$ in $\B(H)$.
\end{proof}

The basic elements we need for \emph{a} proof of the existence of the $\ca$-envelope is Arveson's Extension Theorem and \emph{an exhausting method}, e.g., use of transfinite induction in \cite{DrMc05} (note that throughout the proof we do not use Stinespring's Theorem). In our case the latter is given by the following lemma.

%%%%%%%%%%%%%%%%%%%%%%
\begin{lemma}\label{lemma}
Let $X \subseteq \B(H)$. Then there exists a maximal boundary subsystem for $X$.
\end{lemma}
\begin{proof}
Let $\fI$ be the family of all the boundary subsystems in $\B(H)$. Of course $\fI$ is nonempty since $(0)\in \fI$. An application of Zorn's Lemma will give the maximal element. Indeed, let $\{I_k\}$ be a chain in $\fI$. We set $J=\overline{\cup_k I_k}$ which is a subsystem of $\B(H)$, since $\{I_k\}$ is a chain. We will show that $J$ is a boundary subsystem. Let $\nu \in \bbN$ and $[x_{ij}]\in M_\nu(X)$. Then
\begin{align*}
\nor{[x_{ij}]+M_\nu(J)}
&=\inf\{\nor{[x_{ij}]+[a_{ij}]}\colon a_{ij}\in J\}\\
&= \inf\{\nor{[x_{ij}]+[a_{ij}]}\colon a_{ij}\in \cup_k I_k\}.
\end{align*}
Thus for every $\eps>0$ there are $a_{ij}\in \cup_k I_k$ such that
\[
\nor{[x_{ij}] +M_\nu(J)} \leq \nor{[x_{ij}]+[a_{ij}]} \leq \nor{[x_{ij}]
+M_\nu(J)} +\eps.
\]
For each $(i,j) \in \{1,2,\dots,\nu\}\times \{1,2,\dots,\nu\}$, let $I_{k_{ij}}\in \{I_k\}$ such that $a_{ij}\in I_{k_{ij}}$. Since $\{I_k\}$ is a chain and $I_{k_{ij}}$ are finite in number, there is an $I\in\{I_k\}$ such that $I_{k_{ij}}\subseteq I$, for every $i,j=1,\dots,\nu$. But $I$ is a boundary subsystem, therefore for every $\eps>0$ we get,
\begin{align*}
 \nor{[x_{ij}]}&=\nor{[x_{ij}+I]}=\inf\{\nor{[x_{ij}+b_{ij}]}\colon  b_{ij}\in
 I\|\}\\ &\leq \nor{[x_{ij}+ a_{ij}]} \leq \nor{[x_{ij}]+M_\nu(J)}
 +\eps \leq \nor{[x_{ij}]}+\eps
\end{align*}
Thus $\nor{[x_{ij}]}=\nor{[x_{ij}]+M_\nu(J)}$, so $J$ is an upper bound for the chain $\{I_k\}_k$ and Zorn's Lemma applies.
\end{proof}

By Proposition \ref{pull back proposition} one can easily deduce the following.

%%%%%%%%%%%%%%%%%%%%%%
\begin{corollary}\label{cor contains not}
Let $\I$ be a maximal boundary subsystem for $X$ in $\B(H)$. Then $\B(H)/\I$ contains no non-trivial boundary subsystems for $q_\I(X)$.
\end{corollary}

From now on let us fix a maximal boundary subsystem $\I$ for $X$ in $\B(H)$. Our aim is to prove that $(\B(H)/\I,q_\I)$ is an injective envelope of $X$. In order to do so we use the notion of averages of a ucc map. Let $\phi\colon \B(H) \rightarrow \B(H)$ be a ucc map. For a fixed $k\geq 1$, we define the map
\begin{align*}
m_k(\phi)\colon \B(H) \rightarrow \B(H)\colon  a \mapsto \frac{\phi(a)+\dots+\phi^k(a)}{k}, \foral a\in \B(H).
\end{align*}
Since $\phi$ is linear, $m_k(\phi)$ is also linear. Moreover $m_k(\phi)$ is contractive because
\begin{align*}
\nor{m_k(\phi)(a)}=\nor{\frac{\phi(a)+\dots+\phi^k(a)}{k}} \leq k\nor{\frac{\phi(a)}{k}} = \nor{\phi(a)} \leq \nor{a}.
\end{align*}
For $\nu \geq 1$, we get that $\left(m_k(\phi)\right)_\nu=m_k(\phi_\nu)$. Hence $m_k(\phi)$ is a ucc map for all $k\in \bbN$.

Similarly, if $S$ is an operator system (in $\B(H)$) and $\phi\colon S \rightarrow S$ is a ucc map, then $m_k(\phi)\colon S \rightarrow S$ is a well defined ucc map.

%%%%%%%%%%%%%%%%%%%%%%
\begin{lemma}\label{lemma ext}
Let $X$ be a unital operator space and $S$ be an operator system with $X \subseteq S$. If $\phi\colon S \rightarrow S$ is a ucc map extending $\id_X$ and $h$ is a self-adjoint element in $S$, then the subspace $\Span\{\phi(h)-h\}$ of $S$ is a boundary subsystem for $X$.
\end{lemma}
\begin{proof}
It is immediate that $\Span\{\phi(h)-h\}$ is a self-adjoint subspace of $S$. Fix a $k\in \bbN$ and let $m_k(\phi)$ be the $k$-average of $\phi$. Since $\phi(x)=x$ for all $x\in X$, it is easy to check that $m_k(\phi)(x)=x$ for all $x\in X$. For any $\gl\in \bbC$ we will then have
\begin{align*}
(\dagger)\quad \nor{x + \gl(\phi(h)-h)} &\geq \nor{m_k(\phi)(x + \gl(\phi(h)-h))}\\
& = \nor{x + m_k(\phi)(\gl(\phi(h)-h))}\\
& = \nor{x + \gl m_k(\phi)(\phi(h)-h)}\\
& = \nor{x + \frac{\gl}{k}(\phi^{k+1}(h)-\phi(h))}.
\end{align*}
But,
\begin{align*}
\nor{\frac{\gl}{k}(\phi^{k+1}(h)-\phi(h))} \leq \frac{|\gl|}{k}\left(\nor{\phi^{k+1}(h)} +\nor{h}\right) \leq \frac{2|\gl|}{k}\nor{h} \ \stackrel{k}{\longrightarrow} \ 0.
\end{align*}
Hence, by taking the limit with respect to $k$ to the inequality ($\dagger$) we have
\begin{align*}
\nor{x + \gl(\phi(h)-h)} \geq \nor{x}, \foral x\in X.
\end{align*}
Taking infimum over $\gl\in \bbC$ we deduce that
\begin{align*}
\nor{x + \Span\{\phi(h)-h\}} \geq \nor{x}, \foral x\in X,
\end{align*}
and the last inequality becomes equality by noting that
\begin{align*}
\nor{x} & \geq \inf\{ \nor{x + v}\colon  v\in \Span\{\phi(h)-h\}\}
 = \nor{x + \Span\{\phi(h)-h\}}.
\end{align*}
The same arguments can be repeated for the matrix-norms, by substituting $x$, $m_k(\phi)$ and $\gl(\phi(h)-h)$ with $[x_{ij}]$, $m_k(\phi_\nu)$ and $[\gl_{ij}(\phi(h)-h)]= [\gl_{ij}] \diag\{\phi(h)-h \}$, respectively. Hence $\Span\{\phi(h)-h\}$ is a boundary subsystem for $X$ in $S$.
\end{proof}

%%%%%%%%%%%%%%%%%%%%%%
\begin{lemma}\label{lemma rigidity}
If $\I$ is a maximal boundary subsystem of $X$ in $\B(H)$, then $(\B(H)/\I,q_\I)$ is a rigid extension of $X$.
\end{lemma}
\begin{proof}
By definition $(\B(H)/\I,q_\I)$ is an extension of $X$. Let $\phi\colon \B(H)/\I \rightarrow \B(H)/\I$ be a ucc map such that $\phi(x+\I)=x+\I$ for all $x\in X$. Recall that $\B(H)/\I$ is spanned by its self-adjoint elements. Therefore if $\phi \neq \id$ then there is a self-adjoint element $h\in \B(H)/\I$ such that $\phi(h) \neq h$. Hence by Lemma \ref{lemma ext} the subspace $\Span\{\phi(h)-h\}$ is a non-trivial subsystem for $X$ in $\B(H)/\I$. But this contradicts with Corollary \ref{cor contains not}.
\end{proof}

%%%%%%%%%%%%%%%%%%%%%%
\begin{lemma}\label{injective envelope}
If $\I$ is a maximal boundary subsystem for $X$ in $\B(H)$, then $(\B(H)/\I,q_\I)$ is an injective envelope for $X$.
\end{lemma}
\begin{proof}
We will show that $(\B(H)/\I,q_\I)$ is an injective and rigid extension of $X$. By Lemma \ref{lemma rigidity} it suffices to show that $\B(H)/\I$ is injective.

First, let us fix a ucc map $\gs\colon  \B(H)/\I \rightarrow \B(H)$ that extends the mapping $x+\I\mapsto x$. The existence of $\gs$ is implied by the injectivity of $\B(H)$. Note that $q_\I\circ \gs|_X =\id_X$, thus $q_\I\circ \gs =\id$, by Lemma \ref{lemma rigidity}.

To show that $\B(H)/\I$ is injective, let $Z\subseteq Y$ be operator spaces and $\phi\colon Z \rightarrow \B(H)/\I$ be a ucc map. Then the map $\gs\circ \phi\colon  Z \rightarrow \B(H)$ is a ucc map, thus it extends to a ucc map $\widetilde{\gs\circ \phi}\colon  Y \rightarrow \B(H)$. Then the ucc map $q_\I\circ \widetilde{\gs \circ \phi}\colon  Z \rightarrow \B(H)/\I$ extends $\phi$, since
\begin{align*}
q_\I\circ \widetilde{\gs \circ \phi}(z)= q_\I\circ \gs \circ \phi(z) = \id \circ \phi(z)=\phi(z), \foral z\in Z,
\end{align*}
and the proof is complete.
\end{proof}

%%%%%%%%%%%%%%%%%%%%%%
\begin{remark}
Once the existence of an injective envelope is proved, the existence of the $\ca$-envelope is implied. The key elements are the Choi-Effros' Theorem and Lemma \ref{range in inj env} that will follow. Note that Lemma \ref{injective envelope} gives a concrete picture of the injective envelope which we investigate in Section \ref{interpretation}. Due to this fact, we can replace Choi-Effros' Theorem with a simpler, yet similar, argument (see Remark \ref{simple choi-effros}).
\end{remark}

%%%%%%%%%%%%%%%%%%%%%%
\begin{theorem} \textup{(Choi-Effros)}
Let $S \subseteq \B(H)$ be an injective operator system and let a ucp $\phi\colon  \B(H) \rightarrow S$ be a projection onto $S$. Then setting $a\odot b = \phi(ab)$ defines a multiplication on $S$ and $S$ together with this multiplication and its usual $*$-operation is a $\ca$-algebra.
\end{theorem}

For the next lemma, recall that the \emph{multiplicative domain} of a ucc map $\phi\colon C \rightarrow \B(H)$ is the $\ca$-subalgebra of $C$
\[
C_\phi:= \{a\in C: \phi(a)^*\phi(a)=\phi(a^*a) \text{ and } \phi(a)\phi(a)^*=\phi(aa^*)\}.
\]
The restriction of $\phi$ to $C_\phi$ is a $*$-homomorphism.

%%%%%%%%%%%%%%%%%%%%%%
\begin{lemma}\label{range in inj env}
Let $(S,\iota)$ be an injective envelope for an operator space $X \subseteq \B(H)$. If $\phi\colon \B(H) \rightarrow S$ is a ucc extension of the mapping $x\mapsto \iota(x)$ then $X$ is in the multiplicative domain of $\phi$, with respect to the $\ca$-algebraic structure $(S,\odot)$ induced on $S$ by Choi-Effros' Theorem. Consequently, the restriction of $\phi$ to $\ca(X)$ is a $*$-homomorphism.
\end{lemma}
\begin{proof}
Let $\psi\colon  S \rightarrow \B(H)$ be a ucc extension of the mapping $\iota(x)\mapsto x$. The rigidity of $S$ implies that $\phi \circ \psi =\id$.

By Schwarz inequality we obtain $\psi(\iota(x))^*\psi(\iota(x)) \leq \psi(\iota(x)^*\odot\iota(x))$, hence applying the ucc map $\phi$ we get that
\begin{align*}
\iota(x)^* \odot \iota(x)& = \phi(x)^* \odot \phi(x)  \leq \phi(x^* x)\\
& = \phi(\psi(\iota(x))^* \psi(\iota(x))) \leq \phi\circ \psi(\iota(x)^*\odot \iota(x))\\
&= \id(\iota(x)^* \odot \iota(x)) = \iota(x)^*\odot \iota(x),
\end{align*}
therefore $\phi(x)^*\odot\phi(x)=\phi(x^*x)$, for all $x\in X$. A symmetric calculation shows also that $\phi(x)\odot \phi(x)^*=\phi(xx^*)$, for all $x\in X$, which completes the proof.
\end{proof}

%%%%%%%%%%%%%%%%%%%%%%
\begin{theorem}\label{main theorem}
Let $X\subseteq C=\ca(X)$ be an operator space. Then the \v{S}ilov ideal $I$ exists and $C/I$ is the $\ca$-envelope of $X$.
\end{theorem}
\begin{proof}
Assume that $C \subseteq \B(H)$ and fix a maximal boundary subsystem $\I$ for $X$ in $\B(H)$. By Lemma \ref{range in inj env}, $X$ is in the multiplicative domain of the ucc map $q_\I\colon \B(H) \rightarrow \B(H)/\I$, hence the restriction of $q_\I$ to $C=\ca(X)$ is a $*$-homomorphism. Therefore, $I=\ker(q_\I|_C)$ is an ideal in $C$. Moreover, $I= \ker q_\I \cap C =\I\cap C$, and $C/I\simeq q_\I(C)= q_\I(\ca(X))=\ca(q_\I(X))$.

If $\J$ is a second maximal boundary subsystem for $X$ in $\B(H)$ then $\B(H)/\J$ is also an injective envelope, hence there is a unique ucis and onto map $\Phi\colon \B(H)/\J \rightarrow \B(H)/\I$ that fixes $X$, by Lemma \ref{nec un in en}. Therefore, the following diagram commutes
\begin{align*}
\xymatrix{  & & \B(H)/\J \ar[d]^{\Phi}\\
         C \ar[urr]^{q_\J} \ar[rr]^{q_\I} & & \B(H)/\I
}
\end{align*}
since $\Phi(q_\J(x))=q_\I(x)$ for all $x\in X$, $X$ is in the multiplicative domain of $q_\J$, $q_\J(X)$ is in the multiplicative domain of $\Phi$ and $C=\ca(X)$. Thus $I=\ker q_\I|_{C}= \ker(\Phi\circ q_{\J}|_{C})= \ker q_\J|_{C}= \J \cap C$. Hence, $I= \I \cap C$ for \emph{any} maximal boundary subsystem $\I$ for $X$ in $\B(H)$.

We will show that $I$ is the \v{S}ilov ideal. Since $I \subseteq \I$ then $I$ is also a boundary (ideal) for $X$ in $C$. Let $J$ be a boundary ideal for $X$ in $C$ and let $\J$ be the maximal boundary subsystem that contains $J$. Then $J \subseteq \J \cap C= I$.

The proof is complete by observing that the $\ca$-cover $(q_\I(C),q_\I)$ has the universal property of the $\ca$-envelope. Indeed, let $(B,j)$ be a $\ca$-cover for $X$ with $B=\ca(j(X)) \subseteq \B(K)$. Then by Lemma \ref{range in inj env} the map $j(x)\mapsto q_\I(x)$ extends uniquely to a $*$-epimorphism $\Phi\colon B \rightarrow \ca(q_\I(X))=q_\I(C)$.
\end{proof}

%%%%%%%%%%%%%%%%%%%%%%
\section{The \v{S}ilov ideal}
%%%%%%%%%%%%%%%%%%%%%%

The proof of Theorem \ref{main theorem} gives additional information for the \v{S}ilov ideal which we isolate in the next corollary.

%%%%%%%%%%%%%%%%%%%%%%
\begin{corollary}\label{main corollary}
Let $X \subseteq C=\ca(X) \subseteq \B(H)$. If $\I$ is any maximal boundary subsystem for $X$ in $\B(H)$, then $\I \cap \ca(X)=I$ is the \v{S}ilov ideal.
\end{corollary}

By definition the \v{S}ilov ideal contains the boundary ideals for $X$ in $C$. But this does not ensure that the \v{S}ilov ideal contains also all the boundary subsystems for $X$ in $C$, as it is not obvious that the ideal generated by a boundary subsystem is in turn boundary. Nevertheless this is implied bythe proof of the existence of the \v{S}ilov ideal provided here.

We say that an $a\in (\B(H))_{\text{sa}}$ is a \emph{boundary element} if the operator subsystem $\Span\{a\}$ is boundary for $X$.

%%%%%%%%%%%%%%%%%%%%%%
\begin{corollary}\label{cor contains}
Let $X\subseteq C=\ca(X) \subseteq \B(H)$. Then
\begin{enumerate}
\item The \v{S}ilov ideal contains all the selfadjoint boundary subsystems for $X$ in $C$. Thus it contains the $($closed$)$ linear span of selfadjoint boundary elements.
\item The ideal generated by a boundary subsystem for $X$ in $C$ is also boundary.
\end{enumerate}
\end{corollary}
\begin{proof}
Let $V\subseteq C$ be a selfadjoint boundary subsystem for $X$ in $C$ and $\I$ be a maximal boundary subsystem that contains $V$. Then, by Corollary \ref{main corollary} $V\subseteq \I\cap C =I$, where $I$ is the \v{S}ilov ideal for $X$. Moreover, since $I$ is a boundary ideal we get that the ideal $\sca{V}$, generated by $V$ in $C$, is also a boundary ideal in $I$.
\end{proof}

%%%%%%%%%%%%%%%%%%%%%%
\begin{corollary}
Let $X\subseteq C=\ca(X) \subseteq \B(H)$. Then the \v{S}ilov ideal for $X$ in $C$ is the biggest boundary subsystem for $X$ contained in $C$.
\end{corollary}

Let $X\subseteq \B(H)$ and $\phi\colon \ca(X) \rightarrow \B(H)$ be a ucc map that extends $\id_X$. We define the set of \emph{abnormalities in $\ca(X)$ relevant to $\phi$} as the set $P_\phi=\{\phi(c)-c\in \ca(X) \mid  c\in \ca(X)\}$. It is clear that $P_\phi$ is a selfadjoint linear subspace of $\ca(X)$, since $\phi$ is a ucp map.

For simplicity, we denote the set $\{\phi\colon\ca(X)\rightarrow \B(H) \mid \phi \text{ \textup{extends} } \id_X\}$ by $\ext(\id_X)$.

%%%%%%%%%%%%%%%%%%%%%%
\begin{proposition}\label{union of sets}
Let $X\subseteq \B(H)$. Then the set \ $\bigcup \{ P_\phi \mid  \phi\in \ext(\id)\}$ equals to the \v{S}ilov ideal, and consequently it is closed.
\end{proposition}
\begin{proof}
Let $I$ be the \v{S}ilov ideal of $X$ in $\ca(X)$ and $a\in I$. Then the ideal $\sca{a}$ that is generated by $a$ is in $I$, thus it is a boundary ideal. Hence the restriction of the quotient map $q_{\sca{a}}$ to $X$ is a ucis. Let $\psi\colon C/\sca{a}\rightarrow \B(H)$ be a ucc map that extends the mapping $x+\sca{a}\mapsto x$. Then the mapping $\psi\circ q_{\sca{a}}$ is a ucc map that extends $\id_X$ and $\psi\circ q_{\sca{a}}(a)=0$. Hence $a=a-\psi\circ q_{\sca{a}}(a)\in P_{\psi\circ q_{\sca{a}}}$. Thus $I \subseteq \bigcup \{ P_\phi \mid  \phi\in \ext(\id_X)\}$.

For the converse, let $a\in \cup \{ P_\phi \mid  \phi\in \ext(\id_X)\}$. Then there is a ucc map $\phi\in \ext(\id_X)$ and a $c\in \ca(X)$ such that $a=\phi(c)-c$. For $c_1=\frac{c+c^*}{2}$ and $c_2=\frac{c-c^*}{2i}$, we get that $a= \phi(c_1)-c_1 + i( \phi(c_2)-c_2)\in \sca{\phi(c_1)-c_1} + \sca{\phi(c_2)-c_2}$. By Lemma \ref{lemma ext} the elements $\phi(c_1)-c_1$ and $\phi(c_2)-c_2$ are boundary. Therefore the ideal $\sca{\phi(c_1)-c_1} + \sca{\phi(c_2)-c_2}$ is boundary as the sum of two boundary ideals, by Corollary \ref{cor contains}. Thus $a\in I$.
\end{proof}

%%%%%%%%%%%%%%%%%%%%%%
\begin{remark}
The definition of the sets $P_\phi$ is rather tricky. It refers to elements $c\in \ca(X)$ such that $\phi(c)\in \ca(X)$ for some ucc extension $\phi$ of $\id_X$, and not to elements such that $\phi(c)\in \B(H)$. The reason to be careful is that a ucc extension $\phi$ of $\id_X$ may take values outside $\ca(X)$, even outside an original injective envelope $S$ that contains $X$ (even when $\ca(X)$ is considered as the $\cenv(X)$). On the other hand this assumption seems reasonable as the \v{S}ilov ideal lies in $\ca(X)$. It would be of great interest a result similar to Proposition \ref{union of sets} for sets of the form $\{\phi(c)-c \mid c\in \ca(X)\}$.
\end{remark}

%%%%%%%%%%%%%%%%%%%%%%
\section{Maximal Boundary Subsystems} \label{interpretation}
%%%%%%%%%%%%%%%%%%%%%%

Corollary \ref{main corollary} associates the \v{S}ilov ideal with maximal boundary subsystems. In this section we investigate further these spaces.

%%%%%%%%%%%%%%%%%%%%%%
\begin{proposition}\label{iff proposition}
Let $X \subseteq \B(H)$. Then $\I$ is a maximal boundary subsystem for $X$ in $\B(H)$ if and only if $\I=\ker\phi$ for some ucc map $\phi\colon \B(H) \rightarrow S$, where $(S,\iota)$ is an injective envelope for $X$ and $\phi(x)=\iota(x)$, for all $x\in X$.
\end{proposition}
\begin{proof}
If $\I$ is a maximal boundary subsystem then the appropriate $\phi$ is $q_\I\colon \B(H) \rightarrow \B(H)/\I$. For the converse, let $S$ be an injective envelope and $\phi\colon \B(H) \rightarrow S$ be a ucc map that fixes $X$ elementwise. By Lemma \ref{nec onto in en}, the map $\phi$ is onto $S$. Therefore $\B(H)/\ker\phi \simeq S$ via the ucc map $\widehat \phi\colon  x+\ker\phi \mapsto \phi(x)$. Thus $\ker\phi$ is a boundary operator subsystem (because $\phi$ is also positive) for $X$. Moreover, $\B(H)/\ker\phi$ is also an injective envelope for $X$. Indeed, it suffices to prove that $(\widehat\phi)^{-1}$ is ucc, because then $\widehat \phi$ is ucis. To this end let $\sigma\colon S \rightarrow \B(H)$ be the ucis map that extends the map $\iota(x) \mapsto x$. Then $\widehat \phi \circ q_{\ker\phi} \circ \sigma\colon S \rightarrow S$ and
\begin{align*}
\widehat \phi \circ q_{\ker\phi} \circ \sigma (\iota(x))=\widehat \phi (x+\ker\phi)=\phi(x)=\iota(x), \foral x\in X.
\end{align*}
Therefore by rigidity of $S$ we get that $\widehat \phi \circ q_{\ker\phi} \circ \sigma=\id_S$. Then
\begin{align*}
q_{\ker\phi} \circ \sigma&= \id_{\B(H)/\ker\phi}\circ q_{\ker\phi} \circ \sigma\\
 &= (\widehat\phi)^{-1} \circ \widehat \phi \circ q_{\ker\phi} \circ \sigma = (\widehat \phi)^{-1}\circ \id_S= (\widehat\phi)^{-1}.
\end{align*}
Thus the map $(\widehat \phi)^{-1}$ is ucc, hence $\widehat \phi$ is ucis, therefore $\B(H)/\ker\phi$ is an injective envelope.

If $\I$ is a maximal boundary operator subsystem that contains $\ker\phi$, then $\B(H)/\I$ is also an injective envelope. Lemma \ref{nec un in en} implies then that the canonical onto map
\begin{align*}
\B(H)/\ker\phi \rightarrow \B(H)/\I\colon  a+\ker\phi \mapsto a+ \I,
\end{align*}
is a (necessarily unique) ucis, therefore $\ker\phi=\I$.
\end{proof}

%%%%%%%%%%%%%%%%%%%%%%
\begin{remark}\label{simple choi-effros}
We can always assume that there is an injective envelope $(S,\id)$ for $X \subseteq \B(H)$. (Indeed, if $(S,\iota)$ is an injective envelope and $\phi\colon S \rightarrow \B(H)$ is any extension of the mapping $\iota(x)\mapsto x$, then $\phi(S)$ is an injective envelope containing $\phi(\iota(X))=X$. Note that essentiality of $S$ guarantees that $\phi(S)=\overline{\phi(S)}$.) Let $\I$ be a maximal boundary subsystem for $X$; there is a trivial (but isomorphic) way of inducing a $\ca$-algebraic structure on $\B(H)/\I$ avoiding the use of Choi-Effros' Theorem. We view $\B(H)/\I$ as $\B(H)/\ker\phi$, for a $\phi$ as in Proposition \ref{iff proposition}, and for $a,b\in \B(H)$ we define
\begin{align*}
(a+\ker\phi)\odot (b+\ker\phi): = \phi(a)\phi(b)+ \ker\phi.
\end{align*}
It is well defined and in order to prove the $\ca$-identity, we will use that $\phi$ is a ucp projection on $S$ (by Lemma \ref{nec onto in en}), that $\widehat\phi\colon \B(H)/\ker\phi \rightarrow S$ is a ucis map (by Lemma \ref{nec un in en}), and the $\ca$-identity on $\B(H)$. For $a\in \B(H)$ we get
\begin{align*}
\nor{a+\ker\phi}^2 & =\nor{\phi(a)}^2 = \nor{\phi(a)^*\phi(a)} = \nor{\phi^2(a)^*\phi^2(a)}\\
& \leq \nor{\phi(\phi(a)^*\phi(a))} = \nor{\phi(a)^*\phi(a) +\ker\phi}\\
&= \nor{(a+\ker\phi)\odot (a^*+ \ker\phi )}= \nor{\phi(a)^*\phi(a) +\ker\phi}\\
&= \nor{\phi(\phi(a)^*\phi(a))} \leq \nor{\phi(a)^*\phi(a)}\\
&= \nor{\phi(a)}^2 = \nor{a+\ker\phi}^2.
\end{align*}
\end{remark}

Another way of characterizing maximal boundary subsystems for $X$ in $\B(H)$ is by using \emph{minimal $X$-maps}. A map $\phi\colon \B(H) \rightarrow \B(H)$ is called an \emph{$X$-map} if $\phi$ is ucc and $\phi(x)=x$ for all $x\in X$. Trivially, the kernel of an $X$-map is a boundary subsystem. For example,
\begin{align*}
\nor{x}\geq \nor{x+\ker\phi} \geq \nor{\phi(x)}= \nor{x}, \foral x\in X.
\end{align*}

An $X$-map is called an \emph{$X$-projection}, if it is a projection. We write $\psi \prec \phi$, if $\psi$ is an $X$-projection such that $\psi\circ \phi = \psi = \phi\circ \psi$.

For an $X$-map $\phi$ we can define an \emph{$X$-seminorm} $p_\phi$ on $\B(H)$ such that $p_\phi(a)=\nor{\phi(a)}$, for all $a\in \B(H)$. Unlike in \cite{Pau02}, we write $p_\psi \leq_{c} p_\phi$, if $\psi$ is an $X$-map and $p_{\psi_\nu}([a_{ij}]) \leq p_{\phi_\nu}([a_{ij}])$ for all $a_{ij}\in \B(H)$ and $\nu \in \bbN$.

%%%%%%%%%%%%%%%%%%%%%%
\begin{theorem}
Let $X$ be an operator space in $\B(H)$. Then the following are equivalent
\begin{enumerate}
\item $\I$ is a maximal boundary subsystem for $X$ in $\B(H)$,
\item $\I=\ker\phi$ for some $\phi\colon \B(H) \rightarrow S$, where $(S,\iota)$ is an injective envelope for $X$ and $\phi(x)=\iota(x)$, for all $x\in X$.
\item $\I=\ker\phi$ for some $X$-map $\phi$ such that $p_\phi$ is a $\leq_c$-minimal $X$-seminorm,
\item $\I=\ker\phi$ for some $\prec$-minimal $X$-projection $\phi$ and $(\phi(\B(H)), \phi|_X)$ is an injective envelope for $X$.
\end{enumerate}
\end{theorem}
\begin{proof}
The equivalence $[(1)\Leftrightarrow (2)]$ is Proposition \ref{iff proposition}, and the implications $[(3) \Rightarrow (4)]$ can be derived, for example, by the same arguments as in \cite[proof of Theorem 15.4]{Pau02}. Obviously $[(4) \Rightarrow (2)]$.

For $[(2) \Rightarrow (3)]$, we can assume that $S \subseteq \B(H)$, by Remark \ref{simple choi-effros}. Let $\phi\colon  \B(H) \rightarrow S$ be a map that fixes $X$ pointwise. Hence it is an $X$-projection onto $S$ and $\ker\phi$ is a maximal boundary subsystem for $X$ in $\B(H)$. Let $\psi\colon \B(H) \rightarrow \B(H)$ such that $\nor{\psi_\nu([a_{ij}])}\leq \nor{\phi_\nu([a_{ij}])}$ for all $a_{ij}\in \B(H)$ and $\nu \in \bbN$. Then the mapping $\gs\colon  S \rightarrow \B(H)$, such that $\gs(\phi(a))=\psi(a)$, is a well defined ucc map onto $\psi(\B(H))$ and fixes $X$ elementwise. Moreover, $\ker\psi= \ker\gs\circ \phi \supseteq \ker\phi$. Since $\ker\phi$ is a maximal boundary subsystem and $\ker\psi$ is a boundary subsystem (being the kernel of an $X$-map), this implies that $\ker\phi=\ker\psi$. Therefore the following diagram is commutative
\begin{align*}
\xymatrix{ \B(H)/\ker\phi \ar[rr]^{\widehat\phi} \ar@{=}[d] & & S \ar[d]^{\gs}\\
           \B(H)/\ker\psi \ar[rr]^{\widehat\psi} & & \psi(\B(H))
}
\end{align*}
where the induced quotient maps $\widehat\phi$ and $\widehat\psi$ are ucis. Indeed, this is implied by the fact that $S$ is an injective envelope, hence $\B(H)/\ker\phi=\B(H)/\ker\psi$ is also an injective envelope. Thus, $\gs$ is a ucis map. But then
\begin{align*}
\nor{\psi_\nu([a_{ij}])} =\nor{(\gs\circ\phi)_\nu([a_{ij}])}
=
\nor{\gs_\nu(\phi_\nu([a_{ij}]))}= \nor{\phi_\nu([a_{ij}])},
\end{align*}
for all $a_{ij}\in \B(H)$ and $\nu \in \bbN$; hence $p_\phi$ is $\leq_c$-minimal.
\end{proof}

\bigskip

\noindent {\bf Remark.} We remind that our arguments work for the category of unital operator spaces with ucc maps. If $X$ is an operator system this is the usual category, whereas if $X$ is a non-unital operator algebra, one can prove the analogous results by passing to the unitization of $X$ and by using Meyer's Theorem \cite{Mey01}. However, in general there are results that suggest that the $\ca$-envelope of non-unital operator spaces may not be the natural object (see \cite{Ble01, BlePau01, Zha95}). The author would like to thank M. Anoussis and A. Katavolos for bringing this to his attention.

\bigskip

\noindent {\bf Acknowledgements.}
The author would like to thank MD/PhD-candidate Petros Charalampoudis, Dr.  Christos Chorianopoulos, BSc Konstantinos Mantas and PhD-candidate Vaggelis Skotadis for their comprehensive support in times that we ``fail to follow''. Times when, ``what we can do is to expand the field of consciousness in order to extend the visual field with which we retrieve information from the inner and the outer world''.% (see also \cite{Carl}).

%%%%%%%%%%%%%%%%%%%%%%


\begin{thebibliography}{99}
%%%%%%%%%%%%%%%%%%%%%%

\bibitem{Ar06} W. Arveson, \textit{Notes on the unique extension property}, 2006, http: //math.berkeley.edu/~arveson/Dvi/unExt.pdf.

\bibitem{Ble01} D. P. Blecher, \textit{The Shilov boundary of an operator space and the characterization theorems}, J. Funct. Anal. \textbf{182}(2) (2001), 280--343.

\bibitem{BleLeM04} D. P. Blecher and C. Le Merdy, \textit{Operator algebras and their modules---an operator space approach}, volume 30 of \textit{London Mathematical Society Monographs, New Series}, The Clarendon Press Oxford University Press, Oxford, 2004.

\bibitem{BlePau01} D. P. Blecher, V. I. Paulsen, \textit{Multipliers of operator spaces nad the injective envelope}, Pac. J. Math \textbf{200}(1) (2001), 1--17.

\bibitem{DrMc05} M. A. Dritschel, S. A. McCullough, \textit{Boundary representations for families of representations of operator algebras and spaces}, J. Operator Theory, \textbf{53}(1) (2005), 159--167.

\bibitem{Ham79} M. Hamana, \textit{Injective envelopes of operator systems}, Publ. RIMS Kyoto Univ. \textbf{15} (1979), 773--785.

\bibitem{Mey01} R. Meyer, \textit{Adjoining a unit to an operator algebra}, J. Operator Theory \textbf{46}  (2001), 281--288.

\bibitem{Pau02} V. I. Paulsen, \textit{Completely bounded maps and operator algebras}, volume 78 of \textit{Cambridge Studies in Advanced Mathematics}, The Cambridge University Press, Cambridge, 2002.

\bibitem{Zha95} C. Zhang, \textit{Representations of operator spaces}, J. Operator Theory \textbf{33}(2) (1995), 327--351.


\end{thebibliography}
\end{document}